\documentclass[12pt]{article}
\newcommand{\nsec}[1]{\vspace{-4ex}\section{#1}\vspace{-1ex}}
\newcommand{\diagsm}[1]{\spreaddiagramcolumns{-#1 pc}
\spreaddiagramrows{-#1 pc}}

\usepackage{xypic,amssymb,amsmath,latexsym}

\newcommand{\sitem}{\vspace{-1ex} \item}

\newcommand{\sbib}{\vspace{-1.5ex} \bibitem}
\def\ge{\geqslant}
\def\le{\leqslant}
 \newcommand{\midsq}[1]{\save\go[0,0];[1,1]:(0.5,0)\drop{#1}
   \restore}

\newcommand{\vv}{\mbox{\rule{0.08em}{1.7ex}
\hspace{0.6em}\rule{0.08em}{1.7ex}}}

 \newcommand{\all}{{\mbox{$\rule{0.08em}{1.7ex}\hspace{-0.07em}
\rule[1.5ex]{0.7em}{0.2ex}
\hspace{-0.67em}\rule{0.68em}{0.2ex}\rule{0.08em}{1.7ex}$}}\,}

\newcommand{\Ker}{\mbox{Ker }}

 \newcommand{\labto}[1]{\stackrel{#1}{\longrightarrow} }
\newtheorem{example}{Example}[section]

 \newtheorem{Def}[example]{Definition}

  \newtheorem{prop}[example]{Proposition}

 \newtheorem{thm}[example]{Theorem}

\newtheorem{blank}[example]{\hspace{-0.25em}}

 \newcommand{\sqdiagram}[8]{ \diagram  #1  \rto^{#2} \dto_{#4} &
#3  \dto^{#5} \\ #6    \rto_{#7}  &  #8   \enddiagram }

\newcommand{\sast}{_{\ast}}

 \newcommand{\C}{\v{C}ech}
\newcommand{\B}{\mathsf{B}}
  \newenvironment{proof}{\noindent {\bf Proof} }{ \hfill $\Box$
\mbox{}}

\newcommand{\llabto}[2]{\stackrel{#2}
{\rule[0.54ex]{#1 em}{0.07ex}\hspace{-0.4em}\longrightarrow}}

\begin{document}
\setcounter{page}{187}
 \title{\large \bf  COMPUTING HOMOTOPY TYPES USING
\\ CROSSED  $N$-CUBES OF GROUPS\thanks{This paper (May, 1994; revised,  August, 2006) is an edited
version in {\small \LaTeX\ } of the paper of the same title which
appeared in the {\em Adams  Memorial Symposium on Algebraic
Topology}, Vol 1, edited N. Ray and G  Walker, Cambridge University
Press, 1992, 187-210.  }}

\author{Ronald Brown   \\ University of
Wales, Bangor \\Gwynedd LL57 1UT,  U.K.\\ {\em Dedicated to the
memory of Frank Adams}}
\date{}

\maketitle

\section*{Introduction}
\vspace{-1ex}
The aim of this paper is to explain how, through the work of a
number  of people, some algebraic structures related to
groupoids have yielded  algebraic descriptions of homotopy
$n$-types.  Further, these descriptions are explicit, and in
some cases completely computable, in a way not  possible with
the traditional Postnikov systems, or with other models, such
as simplicial groups.

     These algebraic structures take into account the action of
the  fundamental group.  It follows that the algebra has to be
at least as  complicated as that of groups, and the basic facts
on the use of the  fundamental group in 1-dimensional homotopy
theory are recalled in Section  1.  It is reasonable to suppose
that it is these facts that a higher  dimensional theory should
imitate.

     However, modern methods in homotopy theory have generally
concentrated  on methods as far away from those for the
fundamental group as possible.   Such a concentration has its
limitations, since many problems in the  applications of
homotopy theory require a non-trivial fundamental group  (low
dimensional topology, homology of groups, algebraic $K$-theory,
group  actions, $\ldots$).  We believe that the methods outlined
here continue a  classical tradition of algebraic topology.
Certainly, in this theory non-Abelian groups have a clear role,
and the structures which appear arise  directly from the
geometry, as algebraic structures on sets of homotopy  classes.

     It is interesting that this higher dimensional theory
emerges not  directly from groups, but from groupoids.  In
Sections 2 and 3 we state  some of the main facts about the use
of multiple groupoids in homotopy  theory, including two notions
of {\em higher homotopy groupoid}, and the  related  notions of
{\em crossed complex} and of {\em crossed $n$-cube of groups}.
Theorem \ref{square}  shows how to calculate standard homotopy
invariants of 3-types for the  classifying space of a crossed
square.  We also show in Section 3 how  crossed  $n$-cubes of
groups and the notion of  $n$-cube of fibrations may,  with the
use of the  Generalized Seifert-Van Kampen Theorem due to Brown
and  Loday, 1987a, (Theorem \ref{gvkt}), be used for the
computation of homotopy types in some  practical cases
(Proposition \ref{equiv}).

     An interesting methodological point is that the description
of the  whole  $n$-type has, by these methods, better algebraic
properties than do  the individual invariants (homotopy groups,
Whitehead products, etc.).  As  an application, we give some
explicit results on 3-types, including  computations of
Whitehead products at this level.  In Section 4 we prove a
result from Section 2.  In Section 5, we show that all simply
connected 3-types arise from a crossed square of Abelian groups
(Theorem 5.1).

     Baues, 1989, 1991, also considers algebraic models of
homotopy types  involving non-Abelian groups, and in the second
reference considers  {\em quadratic modules} and {\em quadratic
chain complexes}.  It seems that  the sets of  results of the
two techniques have a non-trivial intersection, but neither  is
contained in the other.  Thus a further comparison, and possibly
integration, of the two types of methods would be of interest.

     Joyal and Tierney have also announced a model of 3-types
using  {\em braided  $2$-groupoids}.  These models are
equivalent to the {\em braided crossed  modules} of  Brown and
Gilbert, 1989, which are there related to simplicial groups and
used to discuss automorphism structures of crossed modules.

\nsec{Groups and homotopy 1-types}      The utility of groups
in homotopy theory arises from the standard  functors
\begin{eqnarray*} &\pi_1  &: (\mbox{spaces with base point}) \to
(\mbox{groups})\\  & \B & : (\mbox{groups})\to   (\mbox{spaces with
base point}) \end{eqnarray*} known as the {\em fundamental group}
and {\em classifying space} functors  respectively.  The classifying
space functor is the composite of geometric realisation and  the
nerve functor  $N$  from groups to simplicial sets.  These functors
have  the properties:
\begin{blank}  There is a natural equivalence of functors
$\pi_1\circ \B  \simeq 1$.\end{blank}   \begin{blank}  If  $G$ is a
group, then  $\pi _i\B G$  is  $0$  if  $i \ne  1$. \end{blank}
  The fundamental group of many spaces may be calculated using
the  Seifert-Van Kampen theorem, or using fibrations of spaces.
Further, if   $X$   is a connected  $CW$-complex and  $G$  is a
group, then there is a natural  bijection  $$[X,\B G] \cong
Hom(\pi_1X,G) , $$ where the square brackets denote pointed homotopy
classes of maps.  It  follows that there is a map  $$X \to
\B\pi_1X$$ inducing an isomorphism of  fundamental groups. It is in
this sense that groups are said to model  homotopy 1-types, and a
computation of a group  $G$  is also regarded as a  computation of
the  1-type of the classifying space  $\B G$.

 A standard block against generalising this theory to higher
dimensions  has been that higher homotopy groups are Abelian. The
algebraic reason for  this is that group objects in the category of
groups are Abelian groups.   This seems to kill the case for a
`higher dimensional group theory', and in  1932 was the reason for
an initial dissatisfaction with \C's definition  of higher homotopy
groups (\C, 1932). Incidentally, \C\  also suggested  the idea of
higher homotopy groups went back to Dehn, who never published  it
(Dieudonn\'e, 1989).  The difficulties of basic homotopy theory are
shown  by the fact that Hurewicz never published the proofs of the
results  announced in his four notes on homotopy groups (Hurewicz,
1935, 1936); that  a proof of the Homotopy Addition Theorem did not
appear in print till Hu,  1953; and that even current proofs of this
basic theorem are not easy (e.g. G.W.Whitehead, 1978).

     It has for some time been established that most of the
theory of the  fundamental group is better expressed in terms of
the fundamental groupoid  (Brown, 1968, 1988) in that theorems:
\begin{itemize}     \sitem  have more natural and convenient
expression;     \sitem  have more powerful statements;
\sitem and have simpler proofs.  \end{itemize}  As an example,
we mention the description in Brown, 1988, of the  fundamental
groupoid of the orbit space of a discontinuous action.  Thus it
is natural to ask: \newline {\em  Can a `better' higher homotopy
theory be obtained using some  notion of `higher homotopy
groupoid'? } \newline  Expectations in this direction were
expressed in Brown, 1967.

     By now, some of the answers to this question seem to be of
the `best  possible situation' kind, and suggest that homotopy
theory is in principle  coincident with a `higher dimensional
group(oid) theory'.  Such a theory is  a significant
generalisation of group theory.  In view of the many
applications of group theory in mathematics and science, the
wider uses of  this generalisation, and the principles
underlying it, need considerable  further study.  For example,
some possibilities are sketched in Brown,  1989b, 1990, Brown,
Gilbert, Hoehnke, and Shrimpton, 1991.  Also, the known
applications in homotopy theory have so far used what seems only
a small  part of the algebra.
\nsec{ Multiple groupoids}
The simplest object to consider as a candidate for `higher
dimensional  groupoid' is an  {\em $n$-fold groupoid}.  This is
defined inductively as a  {\em groupoid object in the category
of $(n - 1)$-fold groupoids}, or  alternatively, as a {\em set
with  $n$  compatible groupoid structures}.  The  compatibility
condition is that if  $\circ _i$   and  $\circ _j$  are two
distinct  groupoid structures, then the {\em interchange law}
holds, namely that  $$(a \circ _i b) \circ _j (c \circ _i d) =
(a \circ _j c) \circ _i (b \circ  _i d)$$  whenever both sides
are defined.  This is often expressed in terms of the  diagram
$$ \left[ \begin{array}{cc} a& b \\ c & d \end{array}  \right]
\diagram  \dto^>>{\textstyle j} \rto & i \\ & \enddiagram $$

Note that Ehresmann, 1963, defines double categories, and the
above  definition is a simple extension of that concept.  The
argument that a  group object in the category of groups is an
Abelian group now yields that  a double groupoid contains a
family of Abelian groups, one for each vertex  of the double
groupoid.  More generally, a basic result is that a double
groupoid contains also two {\em crossed modules over groupoids}
(Brown and  Spencer, 1976).  For example, the {\em horizontal
crossed module} is defined  analogously to the second relative
homotopy group.  It consists in  dimension 2 of the elements of
the form  $$   \def\labelstyle{\textstyle}
\spreaddiagramrows{1.2pc}
\spreaddiagramcolumns{1.2pc}\objectmargin{0pc} \diagram
\dline_{1_V} \rdashed^{\partial m}   \midsq{m} & \dline ^{1_V}
\\ \rline_{1_H} &     \enddiagram  $$   where  $1_H$  and  $1_V$
denote identities for the horizontal and vertical  `edge'
groupoid structures respectively.  In dimension  1  it consists
of the horizontal `edge' groupoid.  The boundary  $\partial m $
of an element $m$  is as shown, and the action is  not hard to
define, as suggested by the following diagram:  $$ m^b =
\def\labelstyle{\textstyle} \spreaddiagramrows{1.2pc}
\spreaddiagramcolumns{1.2pc} \objectmargin{0pc} \diagram
\dline_{1_V} \rdashed^{b^{-1}} \midsq{{\vv} } &
\rdashed^{\partial m} \dline \midsq{m} & \dline \rdashed^b
\midsq{{\vv} }& \dline^{1_V}  \\  \rdashed_{b^{-1}} &
\rline_{1_H} & \rdashed_b & \enddiagram  $$ where $\vv$
denotes a vertical identity.

 The existence of these crossed  modules in any double groupoid,
and the fact that double groupoids can be  constructed from any
given crossed module (Brown and Spencer, 1976),  together
illustrate that double groupoids are in some sense {\em more
non-Abelian than groups}.  This in principle makes them more
satisfactory as  models for two dimensional homotopy theory than
are the second homotopy  groups.  In fact it is known that
crossed modules over groupoids, and hence  also double
groupoids, model all homotopy 2-types (see Mac Lane and
Whitehead, 1950, but note that they use 3-type for what is now
called  2-type).

     One of the features of the use of multiple groupoids is
that they are  most naturally considered as cubical objects of
some kind, since they have  structure in different directions.
Analogous simplicial objects may in  some cases be defined, but
their properties are often difficult to  establish, and are
sometimes obtained by referring to the cubical analogue.   For a
general background to problems on algebraic models of homotopy
types,  see Grothendieck, 1983, although this work does not take
into account the  use of multiple groupoids.

     The first example of which we know of a `higher homotopy
groupoid' was  found in 1974 (see Brown and Higgins, 1978), 42
years after {\C}'s  definition of homotopy groups, namely the
{\em fundamental double groupoid  of a  pair of spaces}.  It is
conveniently expressed in the more general situation  of
filtered spaces as follows (Brown and Higgins, 1981b, as
modified in  Brown and Higgins, 1991, Section 8).  Let
$$X_{\ast}  : X_0 \subseteq  X_1 \subseteq  \cdots  \subseteq
X_n \subseteq   \cdots \subseteq  X $$ be a filtered space.  Let
$RX\sast$  be the {\em singular cubical complex} of   $X\sast$,
consisting for all  $n \ge  0$  of the filtered maps  $I^n\sast
\to   X\sast  $ , where the  standard  $n$-cube  $I^n$  is
filtered by its skeleta, and with the standard  face and
degeneracy maps.  Let  $\rho X\sast$   consist in dimension  $n$
of  the  homotopy classes, through filtered maps and rel
vertices, of such filtered  maps.  (The modification in the 1991
paper is to take the homotopies rel  vertices.)  The standard
gluing of cubes in each direction imposes an extra  structure of
$n$ compositions on  $(RX\sast )_n$  for each  $n \ge  1$.

 It  is a subtle  fact (Brown and Higgins, 1981b) that this
structure is inherited by  $\rho  X\sast$    to give the latter
the structure of  $n$-fold groupoid in each dimension   $n$.
There is also an extra, easily verified structure, on both
$RX\sast $  and   $\rho X\sast $,  namely that of {\em
connections}: these are extra degeneracy operations which  arise
in the cubical context from the monoid structure  max  on the
unit  interval  $I$.  The total structure on  $\rho X\sast $  is
called that of   $\omega $-{\em groupoid}  (Brown and Higgins,
1981a,b).  This gives our first example of a {\em higher
homotopy groupoid}.

     The aim of the introduction of this functor  $$\rho  :
\mbox{(filtered spaces)} \to   (\omega \mbox{-groupoids)} $$ was
that the proof of the usual Seifert-Van Kampen Theorem for the
fundamental group generalised to a corresponding theorem for
$\rho$     (Brown and  Higgins, 1981b).  One main feature of the
proof is that  $\omega $-groupoids,  being  cubical objects, are
appropriate for encoding subdivision methods, since  they easily
allow an `algebraic inverse to subdivision'.  It is not easy to
formulate a corresponding simplicial method.  (See Jones, 1984,
for a  possible approach.)  Another feature crucial in the proof
is the use of the  connections to express facts related to the
Homotopy Addition Lemma.  It  seems that connections are an
important new part of the cubical theory,  since they allow for
`degenerate' elements in which adjacent faces are  identical, as
in the simplicial theory.

     The {\em classifying space}  $\B G$  of an  $\omega
$-groupoid  $G$  is the  geometric  realisation of its
underlying cubical set.  These classifying spaces model  only a
restricted range of homotopy types, namely those which fibre
over a   $K(\pi ,1)$  with fibre a topological Abelian group
(Brown and Higgins, 1991).   Nonetheless, these restricted
models have useful applications.  A principal  reason for this
is the equivalence proved in Brown and Higgins, 1981a,  between
$\omega$-groupoids and the classical tool in homotopy theory of
crossed complex.

     A {\em crossed complex} is a structure which encapsulates
the  properties of  the relative homotopy groups   $
\pi_n(X_n,X_{n-1} ,p) ,~ p \in  X_0 , ~n \ge  2$, for a filtered
space   $X\sast$, together with the boundary maps and the
actions of the fundamental groupoid $\pi_1(X_1,X_0)$ on these
relative homotopy groups.  The notion  was first considered in
the reduced case (i.e. when  $X_0$  is a singleton) by  Blakers,
1948, under the name {\em group system}.  It was studied in the
free  case, and under the name {\em homotopy system}, by
Whitehead, 1949.  The term  {\em crossed complex} is due to
Huebschmann, 1980, who used crossed  $n$-fold  extensions to
represent the elements of the  $(n + 1)$st cohomology group of
a group (see also Holt, 1979, Mac Lane, 1979, Lue, 1981), and to
determine  differentials in the Lyndon-Hochschild-Serre spectral
sequence  (Huebschmann, 1981).  Lue, 1981, gives a good
background to the general  algebraic setting of crossed
complexes.  Crossed complexes have the  advantage of being able
to include information on chain complexes with a  group  $G$  of
operators and on presentations of the group  $G$.  The  category
of crossed complexes also has a monoidal closed structure (Brown
and  Higgins, 1987), which is convenient for expressing
homotopies and higher  homotopies.

     The Generalized Seifert-Van Kampen Theorem for the
fundamental   $\omega $-groupoid of a filtered space (Brown and
Higgins, 1981b) implies  immediately a similar theorem for the
fundamental crossed complex, and this  theorem has a number of
applications, including the Relative Hurewicz  Theorem.  The
latter theorem is thus seen in a wider context, related to
excision, and in a formulation dealing initially with the
natural map   $\pi _n(X,A) \to   \pi _n (X \cup  CA)$.  This
formulation was a model for  the $(n + 1)$-ad  Hurewicz theorem
(Brown and Loday, 1987b).  Other recent applications of  crossed
complexes are given in Baues, 1988, 1991, Brown and Higgins,
1987,  1989, 1991, Baues and Brown, 1990, Baues and Conduch\'e,
1991.

     More general algebraic models related to groupoids are
associated not  with filtered objects but with  $n$-cubes of
objects.  Let  $\langle n  \rangle$   denote the  set  $\{
1,2,\ldots ,n\}$.  An  $n$-{\em cube}  $C$  in a category  $\cal
C$  is a commutative  diagram with vertices  $C_A$  for  $A
\subseteq  \langle n \rangle$   and  morphisms $C_A \to C_{A\cup
\{ i \} } $  for  $A \subseteq  \langle n \rangle  , ~i \in
\langle n \rangle$, and  $i  \not\in  A$. In particular, a
1-cube is a morphism, and a  2-cube is a commutative square.

     Let  $X_{\ast}$   be an  $n$-cube of pointed spaces.
Loday, 1982,  defines the  {\em fundamental cat$^n$-group}  $\Pi
X_{\ast}$.  (We are following the  terminology and  notation of
Brown and Loday, 1987a.)  Here, a {\em cat$^n$-group} may be
defined to  be an  $n$-fold groupoid in the category of groups.
Alternatively, it is an   $(n + 1)$-fold groupoid in which one
of the structures is a group.  (This is  one of several
equivalent definitions considered in Loday, 1982.)

     For simplicity, we describe  $\Pi X_{\ast}$   in a special
case, namely  when  $X_{\ast}$   arises from a pointed  $(n +
1)$-ad  ${\cal X} = (X;X_1,\ldots ,X_n)$  by the  rule:  $X_{
\langle n \rangle}  = X$  and for  $A$  properly contained in
$\langle n \rangle$,   $X_A = \bigcap _{i \not\in A} X_i$, with
maps  the inclusions.  Let  $\Phi$   be the space of maps  $I^n
\to   X$  which take  the  faces of  $I^n$  in the  $i$th
direction into  $X_i$.  Notice that  $\Phi$  has the  structure
of $n$ compositions derived from the gluing of cubes in each
direction.  Let  $\ast  \in  \Phi$   be the constant map at the
base point.   Then  $G = \pi_1(\Phi ,\ast )$ is certainly a
group.  Gilbert, 1988, identifies   $G$ with  Loday's  $\Pi
X_{\ast}$, so that Loday's results, obtained by methods of
simplicial  spaces, show that  $G$  becomes a cat$^n$-group.  It
may also be shown that  the  extra groupoid structures are
inherited from the compositions on  $\Phi$.   It  is this
cat$^n$-group which is written  $\Pi {\cal X}$   and is called
the  {\em fundamental  cat$^n$-group of the  $(n + 1)$-ad ${\cal
X}$}.  This construction of  Loday is  our second example of a
{\em higher homotopy groupoid}.  We emphasise that the
existence of this structure is itself a non-trivial fact,
containing  homotopy theoretic information.  Also the results of
Gilbert, 1988, are for  the case of  $n$-cubes of spaces.

     The nerve  $NG$  mentioned in Section 1 may be defined, not
only for a  group but also for a groupoid  $G$, to be in
dimension  $i$  the set of  groupoid maps  $\pi_1(\Delta
^i,\Delta^i_0) \to   G$ , where  $\Delta^i_0$   is the set of
vertices of the  $i$-simplex  $\Delta^i$ .  It follows by
iteration that  $N$  defines also a  functor
$$((n+1)\mbox{-fold groupoids}) \to   ((n+1)\mbox{-simplicial
sets}).$$  Hence there is a {\em classifying space functor}  $$\B :
(\mbox{cat$^n$-groups}) \to   (\mbox{pointed spaces}) $$ defined as
the composite of geometric realisation and the nerve functor to $(n
+ 1)$-simplicial sets.  Loday, 1982, proves that if  $G$  is a
cat$^n$-group,  then  $\B G$  is  $(n + 1)${\em -coconnected}, i.e.
$\pi _iBG = 0$  for  $i >  n + 1$. He also shows, with a correction
due to Steiner, 1986, that if $X$  is a  connected,  $(n +
1)$-coconnected  $CW$-complex, then there is a  cat$^n$-group  $G$
such that  $X$  is of the homotopy type of  $\B G$.  In fact, Loday
gives an  equivalence between a localisation of the category of
cat$^n$-groups and the pointed homotopy category of connected, $(n +
1)$-coconnected $CW$-complexes.   This can be put provocatively as
\begin{center} $(n + 1)$-fold groupoids model all homotopy $(n +
1)$-types.  \end{center}  That is, the generalisation from
groups or groupoids to  $(n + 1)$-fold  groupoids is as good for
modelling homotopy types as might be expected.   This result
also shows the surprising richness of the algebraic structure
of  $(n + 1)$-fold groupoids.

     There is an important structure equivalent to that of
cat$^n$-groups,  namely that of {\em crossed  $n$-cubes of
groups} (Ellis and Steiner, 1987).   The  main intuitive idea is
that a crossed  $n$-cube of groups is a crossed module  in the
category of crossed  $(n - 1)$-cubes of groups.  This leads to
the  following definition ({\em loc. cit.}).   \begin{Def}{\em
Let $\langle n \rangle$   denote the set   ${1,2,\ldots ,n}$.  A
{\em crossed $n$-cube  of groups} is a family of groups,  $M_A ,
A \subseteq  \langle n \rangle$, together with morphisms  $\mu
_i : M_A \to   M_{A\setminus \{ i \} }  , ~(i \in  \langle n
\rangle   , ~A \subseteq \langle n \rangle  )$,  and functions
$h : M_A \times  M_B \to   M_{A \cup B}  , (A,B \subseteq
\langle n  \rangle )$,  such that if  $^ab$  denotes  $h(a,b)b$
for  $a \in  M_A$  and  $b \in   M_B$  with $A \subseteq  B$,
then for  $a,a' \in  M_A , ~b,b' \in  M_B , c \in  M_C$  and
$i,j \in   \langle n \rangle$, the following hold:
\begin{enumerate}         \sitem  $ \mu _ia = a$  if  $i \not\in
 A$         \sitem $ \mu _i\mu _ja = \mu _j\mu _ia $
\sitem $ \mu _ih(a,b) = h(\mu _ia,\mu _ib)$          \sitem $
h(a,b) = h(\mu _ia,b) = h(a,\mu _ib) $ if $ i \in  A $ and $ i
\in  B $     \sitem $  h(a,a') = [a,a'] $        \sitem $
h(a,b) = h(b,a)^{-1}  $        \sitem $  h(a,b) = 1$  if $ a = 1
$ or $  b = 1 $        \sitem $  h(aa',b) = ^ah(a',b)h(a,b)
   $       \sitem $  h(a,bb') = h(a,b) \;^bh(a,b')         $
   \sitem $ ^ah(h(a^{-1},b),c) \;^ch(h(c^{-1},a),b) \;^bh(h(b
^{-1} ,c),a) = 1 $        \sitem $ ^ah(b,c) = h(^ab,^ac)$  if $
A \subseteq  B  \cap  C .  $  \end{enumerate}  }\end{Def}

     A {\em morphism} of crossed  $n$-cubes  $(M_A) \to   (N_A)$
is a  family of morphisms  of groups  $f_A : M_A \to   N_A  (A
\subseteq  \langle n \rangle )$  which  commute with the maps $\mu
_i$ and  the functions  $h$.  This gives us a category $Crs^ngp$.
Ellis and Steiner,  1987, show that this category is equivalent to
that of cat$^n$-groups, and  this is the reason for the choice of
structure and axioms in Definition  2.1.  This equivalence shows
that there is a {\em classifying space functor}  $$\B : Crs^ngp \to
Top . $$ This functor would be difficult to describe directly.  (See
Porter, 1993,  for a different account of such a functor.)  The
results for cat$^n$-groups  imply that a localisation of the
category $Crs^ngp$  is equivalent to the  homotopy category of
pointed, connected,  $(n + 1)$-coconnected  $CW$-complexes.

    The
{\em fundamental crossed $n$-cube of groups functor}  $\Pi '$
is  defined from  $n$-cubes of pointed spaces to crossed
$n$-cubes of groups:  $\Pi  'X_{\ast}$   is simply  the crossed
$n$-cube of groups equivalent to the cat$^n$-group  $\Pi
X_{\ast}$. It is  easier to identify  $\Pi '$  in classical
terms in the case  $X_{\ast}$   is  the $n$-cube  of spaces
arising as above from a pointed  $(n + 1)$-ad  ${\cal X}  =
(X;X_1,\ldots ,X_n)$. That is, let  $X_{ \langle n \rangle } =
X$  and for  $A$  properly contained in  $\langle n \rangle$
let   $X_A = \bigcap _{i \not\in A}  X_i$. Then $M = \Pi '{\cal
X}$   is given as follows (Ellis and Steiner,  1987):
$M_{\emptyset}  = \pi_1(X )$ ; if  $A = {i_1,\ldots ,i_r}$, then
$M$   is the homotopy  $(r + 1)$-ad group   $\pi _{r+1}(X ;X_1
\cap  X_i ,\ldots ,X_n   \cap  X_i )$;  the maps  $\mu$    are
given by the usual boundary maps;  the  $h$-functions are  given
by  generalised Whitehead products.  Note that whereas these
separate elements of structure had all been considered
previously, the aim  of this theory is to consider the whole
structure, despite its apparent  complications.  This global
approach is necessary for the Generalized  Seifert-Van Kampen
Theorem, stated below.  That $\Pi ' {\cal X}$   satisfies  the
laws for  a crossed  $n$-cube of groups follows immediately
since  $\Pi '{\cal X}$   is  the crossed  $n$-cube of groups
derived from the cat -group $ \Pi X_{\ast}$.  From now  on, we
abbreviate  $\Pi '$  to  $\Pi$  , the meaning being clear from
the context.

  A crossed  $n$-cube of groups  $M$  gives rise to an  $n$-cube
of crossed   $n$-cubes of groups  $\all M $ where   $$((\all
M)(A))_B  = \left\{ \begin{array}{ll}M_B & \mbox{if }
A'\subseteq B \\1 & \mbox{otherwise} \end{array} \right.  $$
      Then  $\B\all M$  is an  $n$-cube of spaces.  The
generalisation to this context of  the result on the fundamental
group of the classifying space of a group is  that there is a
natural isomorphism of crossed  $n$-cubes of groups  $$\Pi
B{\all M} \cong   M . $$ (See Loday, 1982, for the cat$^n$-group
case, and Brown and Higgins, 1981b,  1991, for the analogous
crossed complex case.)  This result confirms the  appropriate
nature of the axioms (1)-(11) of Definition 2.1.

A description of the homotopy groups of  $\B G$   for a
cat$^n$-group $G$   has  been given in Loday, 1982, in terms of the
homology groups of a non-Abelian  chain complex. This, with some
extra work, yields a result on the homotopy  invariants of the
classifying space of a crossed square (i.e. a crossed 2-cube of
groups).  It is useful first to give the axioms for this in a
different notation.

     A crossed square (Loday, 1982) consists of a commutative
square of  morphisms of groups
\begin{equation}{\sqdiagram{L}{\lambda}{M}{\lambda
'}{\mu}{N}{\nu}{P}}\tag{2.2} \end{equation}   together with actions of
$P$  on the groups  $L,M,N$,  and a function   $h : M \times  N
\to   L $.  This structure shall satisfy the following  axioms,
in which we assume that  $M$  and  $N$  act on  $L,M,N$  via  $P
$:

\vspace{1ex}
\noindent  (2.3)(i) the morphisms  $\lambda ,\lambda ',\mu ,\nu$
and  $\mu \lambda  =  \nu \lambda ' $ are crossed modules and
$\lambda$  and  $\lambda '$  are  P-equivariant;

\noindent (ii)$  h(mm',n) = h( m', n)h(m,n) , h(m,nn') =
h(m,n)h( m, n') ;$

\noindent  (iii) $\lambda h(m,n) = m m   , \lambda 'h(m,n) =  nn
; $

\noindent  (iv) $h(\lambda l,n) = l l   , h(m,\lambda 'l) =  ll
; $

\noindent  (v)  $ h( m, n) =  h(m,n) ; $

\noindent  for all  $l \in  L ,
m,m' \in  M , n,n' \in  N , p \in  P .$
\vspace{1ex}

     We now describe the homotopy groups of $\B G$   for a crossed
square  $G$   as above.  The first part of the following result is a
special case of  results in Loday, 1982. \addtocounter{example}{2}
\begin{thm} \label{square}  Let  $G$  be the crossed square (2.2).
Then the homotopy  groups of  $\B G$   may be computed as the
homology groups of the non-Abelian  chain complex \begin{equation}
L\llabto{2}{(\lambda^{-1},\lambda ')}  M \rtimes  N  \llabto{1}{\mu
\ast \nu}L P  \qquad \qquad \tag{2.5} \end{equation}   where  $\mu
\ast \nu  : (m,n) \mapsto (\mu m)(\nu  n) .$ This implies that
\begin{equation}\pi _i\B G \cong
\begin{cases}P/(\mu M)(\nu N) & \mbox{if }  i =1 \\
 (M\times _P N)/\{ (\lambda l, \lambda ' l) : l \in L \} &
\mbox{if } i = 2 \\
( \Ker \lambda ) \cap ( \Ker  \lambda ') &
\mbox{if }  i = 3 \\
 0 & \mbox{if }  i \ge 4.  \end{cases}
        \tag{2.6} \end{equation}
  Further, under these isomorphisms, the composition  $\eta ^* :
\pi_2 \B G \to   \pi_3 \B G $ with  the Hopf map  $\eta  : S^3 \to
S^2 $ is induced by the function $ M \times  _P N \to  L ,~ (m,n)
\mapsto h(m,n)$,  and the Whitehead product  $\pi_2 \times  \pi_2
\to \pi_3$ on $\B G$  is  induced by the function $((m,n),(m',n'))
\mapsto h(m',n)h(m,n')$.  The first  Postnikov invariant of  $\B G$
is the cohomology class determined by the crossed module  $$ (M
\rtimes N)/Im(\lambda ^{-1} ,\lambda ')  \llabto{1}{\mu \ast \nu} P.
$$
\end{thm}       We will explain the proof of this result in
Section 4. \nsec{ $n$-cubes of fibrations}      As in Brown and
Loday, 1987a, an  $n$-cube of maps  $X_{\ast}$   yields an
$n$-cube  of fibrations  $\overline{ X}\sast $.  (See Edwards and
Hastings, 1976, Cordier and  Porter,  1990.)  Following Steiner,
1986, we parametrize this as a commutative  diagram consisting of
spaces  $X_{A,B}$  ($A,B$  disjoint subsets of $\langle n
\rangle$  ) and fibration sequences
 \begin{equation}
\overline{X}_{A\cup \{ i\} ,B} \to   \overline{X}_{A,B} \to
\overline{X}_{A,B\cup \{ i\} } , A  \cap  B =  \emptyset , i \in
\langle n \rangle \setminus  (A \cup  B) \tag{3.1}\end{equation}
 \addtocounter{example}{1}  The  $n$-cube
of  fibrations  $(\overline{X}_{A,B})$ contains an  $n$-cube of
spaces   $\overline{X}_{\emptyset ,*}$ homotopy equivalent to
$X\sast$  (i.e. there is a morphism $X_{\ast}  \to
\overline{X}_{\emptyset ,*}$  consisting  of homotopy equivalences
$X_{B} \to   \overline{X}_{\emptyset ,B} )$.  The  $n$-cube of
maps  $X_{\ast}$   is called {\em connected} if all the spaces
$\overline{X}_{A,B}$  are path-connected.

     Just as the Seifert-Van Kampen Theorem enables one to
compute the  fundamental group of a union of connected spaces,
so the Generalised  Seifert-Van Kampen Theorem (GSVKT) enables
one to compute the fundamental  crossed  $n$-cube of a union of
connected  $n$-cubes.  This result is Theorem  5.4 of Brown and
Loday, 1987a, where it is proved by induction on $n$  starting
with the usual SVKT.  It may be restated in terms of crossed
$n$-cubes of groups, rather than cat$^n$-groups, as follows.
\begin{thm}  \label{gvkt}  Let  $X_{\ast}$   be an  $n$-cube of
spaces, and suppose that   ${\cal U} = \{ U^{\lambda} \}$  is an
open cover of the space $X _{<n>}$ , such that  ${\cal U}$  is
closed under  finite  intersections.  Let $ {\cal U} ^{\lambda}
$   be the  $n$-cube of spaces obtained from  $X_{\ast}$   by
inverse images of the $ U^{\lambda}$.  Suppose that each  ${\cal
U} ^{\lambda}  $  is a connected  $n$-cube  of spaces.  Then:

\noindent {\rm (C)}: the $n$-cube  $X_{\ast}$   is connected,
and

\noindent {\rm (I)}: the natural morphism of crossed  $n$-cubes
of groups $$ \mbox{colim}^{\lambda}~ \Pi {\cal U} ^{\lambda}
\to    \Pi X_{\ast}$$  is an isomorphism.  \end{thm}      The
colimit in this theorem is taken in the category of crossed
$n$-cubes of groups, and so the validity of (I) confirms again
that the  axioms for crossed  $n$-cubes of groups are well
chosen.

The connectivity statement (C) of this theorem generalises the
famous  $(n + 1)$-ad connectivity theorem, which is usually
regarded as a difficult  result (at the time of writing, no
recent proof is in print except that  referred to here).  Of
course, the connectivity result is related to the  fact that a
colimit of zero objects is zero.

   The isomorphism statement (I) implies the characterisation by
a  universal property of the critical group of certain $ (n +
1)$-ads.  (See  Brown and Loday, 1987b, for the general
procedure and explicit results on  the triad case, using a
non-Abelian tensor product, and Ellis and Steiner,  1987, for
the general case.)  The earlier result in this area is in
Barratt  and Whitehead, 1952, but there the assumption is made
of simply connected  intersection, and the result is proved by
homological methods, so that it  has no possibility for dealing
with the occurrence of a non-Abelian   $(r + 1)$-ad homotopy
group.  It is clearly advantageous to see the Barratt  and
Whitehead result, including the $(n + 1)$-ad connectivity
theorem, as a  special case of a theorem which has other
consequences, for example an   $(n + 1)$-ad Hurewicz theorem
(Brown and Loday, 1987b).

     These results, with Theorem \ref{square}, illustrate how
situations in homotopy  theory may require constructions on
non-Abelian groups for the convenient  statement of a theorem,
let alone its proof.  The methods of crossed   $n$-cubes of
groups give a (largely unstudied) range of new constructions in
group theory.

     Theorem \ref{gvkt} allows in some cases for the computation
of the fundamental crossed  $n$-cube of groups  $\Pi X_{\ast}$
of an $n$-cube of  spaces $X_{\ast}$. We  now consider to what
extent it also  allows  computation  of  the   $(n + 1)$-type of
the space  $X _{<n>}$.

     Let  $X_{\ast}$   be a connected  $n$-cube of spaces, and
let   $X = X _{<n>}$. It is  proved in Loday, 1982, that there
is an  $n$-cube  of  fibrations   $Z\sast$  and maps of
$n$-cubes of fibrations

                  $$  \overline{X}  \stackrel{f}{\longleftarrow}
Z\sast   \labto{g}  \overline{  \B\all (\Pi X\sast )}$$  such that
$f$  is a level weak homotopy equivalence and  $g$ induces an
isomorphism of  $\pi_1$  at each level.  Assume now that  $X$  is of
the homotopy  type of a  $CW$-complex.  Then from  $f$  and  $g$ we
obtain a map   $$\phi  : X \to   \B\Pi X_{\ast}$$   inducing an
isomorphism of $ \pi_1$, namely the composite, {\em at  this level},
of  $g$  with a homotopy inverse of  $f$, and with the map  $X_{<n>}
\to \overline{X}_{<n>}$. We do not  expect  $\phi$   to be a
homotopy equivalence in general, since  the  $n$-cube of fibrations
$\overline{\B\all (\Pi X_{\ast} )}$  has special properties not
necessarily  satisfied by $\overline{X}\sast$.

 We say an  $n$-cube  of  spaces   $X_{\ast}$    is  an  {\em
Eilenberg-Mac Lane} $n$-cube of  spaces if it is  connected  and all
the  spaces   $\overline{X}  _{A,\empty}$  are spaces of type $K(\pi
,1)$.  A chief example of this is the $n$-cube of spaces $\B\all  M$
derived from a crossed  $n$-cube of groups. In   fact, $(\B\all
M)_{A,B}$  is not  only  path-connected but  also   $(|B| +
1)$-coconnected.  This  $n$-cube of fibrations may also be
constructed directly  in terms of the structure of  $M$, using the
techniques of Loday, 1982.       We have the following result.
\begin{prop} \label{equiv}  Let $X_{\ast}$   be a connected
$n$-cube of  spaces such that $X _{<n>}$  is of the homotopy type of
a  $CW$-complex.  Suppose that for  $A,B \subseteq   \langle n
\rangle$,  such that  $A \cap  B = \emptyset  ,~ i \in   \langle  n
\rangle  \setminus ( A \cup  B)$, and  $r = |B|$, the induced
morphism  $\pi _{r+2}\overline{X}_{A,B} \to   \pi
_{r+2}\overline{X}_{A,B\cup \{ i\} }$  is  zero.  Then the canonical
(up to  homotopy)  map
 $\phi   :  X  _{<n>}  \to    \B\Pi  X_{\ast}$   is an  $(n +
1)$-equivalence. \end{prop}  \begin{proof}   This is a simple
consequence of the five lemma applied by  induction  on  $|B|$ to
the maps of homotopy exact sequences of the fibration sequences
(3.1) of the  $n$-cubes of fibrations $\overline{X}\sast$  and
$\overline{\B\all           (\Pi
   {\cal           X}          )}$.               \end{proof}
\begin{example}{\em   Let  $M$ and  $N$  be normal subgroups of
a  group  $P$, and let  the space  $X$  be given as the homotopy
pushout
$$\diagsm{1.2} \sqdiagram{K(P,1)}{}{K(P/M,1)}{}{}{K(P/N,1)}{}{X} $$ Brown
and Loday, 1987a, apply the case  $n = 2$  of Theorem \ref{gvkt}
to show that  the above square of spaces has fundamental crossed
square given by the  `universal' crossed square
\begin{equation}   { \diagsm{1.2} \sqdiagram{M\otimes
N}{}{M}{}{}{N}{}{P}} \label{tens}\tag{3.5} \end{equation}  where  $M
\otimes  N$  is the non-Abelian tensor product ({\em loc. cit.}),
with generators  $m \otimes n$ for  $m \in  M$  and  $n \in  N$  and
relations   \begin{align*} mm' \otimes  n &= (^mm' \otimes{}  ^mn)(m
\otimes  n) ,\\ m \otimes  nn' &= (m \otimes  n)(^nm \otimes
{}^nn')\end{align*}   for all  $m,m' \in  M , n,n' \in  N$.  The
$h$-map of this crossed square is   $(m,n) \mapsto m \otimes  n$. It
follows from Proposition 3.3 that  the 3-type of  $X$ is  also given
by this crossed square.  This result has been stated in Brown,
1989b, 1990, and we have now given the proof.  Note that Theorem
\ref{square} allows  one to compute  $\eta   : \pi_2  \to \pi_3$ and
the Whitehead product map   $\pi_2  \times  \pi_2 \to   \pi_3 .$
}\end{example}

    By contrast, the Postnikov description of the 3-type of  $X$
requires  the description of the first  $k$-invariant
$$k^{(3)} \in  H^3(P/MN,(M  \cap  N)/[M,N]) ,$$  which in this
case is represented by the crossed module  $M \circ   N \to
P$, where  $M \circ  N$  is the coproduct of the crossed
$P$-modules  $M$   and  $N$  (see Brown,  1984, and also Gilbert
and Higgins, 1989).  This  $k$-invariant determines (up to
homotopy) a space  $X^{(2)}$, which could be taken to be the
classifying  space of the above crossed module, constructed
either by regarding the  crossed module as a crossed 1-cube of
groups, or as in Brown and Higgins,  1991.  One then needs a
second Postnikov invariant  $$k^{(4)} \in  H^4(X^{(2)},Ker(M
\otimes  N \to   P)) .$$  This description of the 3-type of  $X$
is less explicit than that given by  the crossed square
(\ref{tens}), from which we obtained the homotopy groups and
the action of  $\pi_1$  in the first place.  Note  also  that
if   $M , N , P$ are  finite, then so also is  $M \otimes  N$
(Ellis, 1987), so that in  this case the crossed square
(\ref{tens}) is finite.

  As an example, in this way one finds that if  $P = M = N$  is
the  dihedral group  $D_n$  of order  $2n$, with generators  $x$
 and   $y$  and  relations  $x^2 = y^n = xyxy = 1$, where $n$ is
even, then the suspension   $SK(D_n,1)$   of   $K(D_{n},1)$
has   $\pi_3$    isomorphic  to  $({\mathbb Z}_{2})^4$
generated by the elements  of   $D_n  \otimes   D_n$: $$   x
\otimes  x , (x \otimes  y)^{n/2} , y \otimes  y , (x \otimes
y)(y \otimes  x).$$   Further,  $\eta ^ *(\overline{x}) = x
\otimes  x , \eta ^* (\overline{y}) = y \otimes   y$, where
$\overline{x}$ and  $\overline{y}$ denote the  corresponding
generators of  $\pi_2 SK(D_{n},1)  =  (D_{n})^{{\rm  ab}}$  (if
$n$ is odd, only  the $x \otimes  x$  term appears in  $\pi_3$).
 The element  $(x \otimes  y)(y \otimes  x)$ is the  Whitehead
product  $[\overline{x}, \overline{y}]$.  Other  computations
of  $\eta^ * $  and of Whitehead  products at this level in
spaces  $SK(P,1)$  may be deduced from the  calculations of
non-Abelian tensor products given in Brown, Johnson and
Robertson, 1987. (This paper covers the case of dihedral,
quaternionic,  metacyclic and symmetric groups, and all groups
of  order   $\le   31$.)       Problems in this area are given
in Brown, 1990.
\nsec{Proof of Theorem 2.4}      We now
explain the results on  $\eta^ *$  and Whitehead products in the
 second part of Theorem \ref{square}.       Let  $G$  be the
crossed square (2.2).  Then the square of crossed  squares
$\all G$  may be written in abbreviated form as follows:
\begin{equation}  \begin{array}{ccc}  {\left( \begin{array}{cc}1
& 1  \\1 &P    \end{array} \right) }& \longrightarrow & {\left(
\begin{array}{cc}1 & 1  \\N  & P    \end{array} \right) } \\
\downarrow & & \downarrow \\ {\left( \begin{array}{cc}1 & M  \\1
&P    \end{array} \right) } & \longrightarrow &    {\left(
\begin{array}{cc}L & M  \\N &P    \end{array} \right)
}\end{array} \tag{4.1} \end{equation} Let us write  $Y\sast $  for
the square of spaces  $\B\all G$.  Then        $\Pi Y$  is
isomorphic to the original crossed square  $G$ .  Further the 2-cube
of fibrations   $\overline{Y} \sast$   associated to  $Y\sast $  is
homotopy equivalent to the following diagram:  \begin{equation}
{\diagsm{1.1} \diagram BL \rto \dto & BM \rto \dto & B(L \to M) \dto  \\
   BM \rto \dto & BP \rto \dto & B(M \to P)\dto  \\  B(L\to N)
\rto & B(M \to P)\rto  &  B(G) \enddiagram   }\tag{4.3} \end{equation}
For a general square of spaces  $X_{\ast}$   as follows
\begin{equation} {  \diagsm{1.1}
\sqdiagram{C}{f}{A}{g}{a}{B}{b}{X} }\tag{4.4} \end{equation}
the associated 2-cube of fibrations is equivalent to the
following diagram  \begin{equation}{\diagsm{1.1} \diagram F( X \sast ) \rto
\dto & F(g) \rto \dto  &F(a) \dto \\ F(f)  \dto \rto & C \rto
\dto & A \dto \\ F(b) \rto & B \rto &X \enddiagram }\tag{4.5}
\end{equation} where each row and column is a fibration sequence.
 So we deduce  the second part of Theorem
\ref{square}  from the  following more general result. \addtocounter{example}{5}
\begin{prop}  Let  $X_{\ast}$   be the square of pointed  spaces
as in (4.4).   Suppose that the induced morphism
$\pi_2 C \to   \pi_2 X$  is  zero.  Then there is a commutative
diagram  \begin{equation} {\diagsm{1.1} \diagram \pi_2 X \dto_{\eta^*} &&
\llto  _{\delta  '}  \dto ^{h'}  \pi_1F(f) \times _{\pi _1 C}
\pi _1 F(g) \\  \pi_3  X  \rto_{\partial} & \pi_2 F(a)
\rto_{\partial '} & \pi_1 F({\bf X})  \enddiagram }\tag{4.7}\end{equation}
in which  $\delta '$  is defined by a  difference construction,
$\partial  , \partial '$  are boundaries in homotopy  exact
sequences  of fibrations,  $\eta ^*$  is induced by composition
with the  Hopf map  $\eta$,  and $h'$  is the restriction  of
the   $h$-map  of  the  crossed  square  $\Pi X_{\ast}$.
\end{prop}  \begin{proof}   This result is a refinement of Lemma
4.2 of Brown and Loday,  1987a.  It is proved by similar
methods.  One first considers the suspension square  of  $S^1$:
$$\diagsm{1.2} \sqdiagram{S^1}{}{E^2_+}{}{}{E^2_{\mbox{-}}}{}{S^2} $$
The fundamental crossed square of this suspension square is
determined by  Theorem 3.2, compare Example 3.4, as in Brown and
Loday, 1987a, and is  $$\diagsm{1.1}\sqdiagram{{\mathbb Z}}{0}{{\mathbb
Z}}{0}{1}{{\mathbb Z}}{1}{{\mathbb Z}} $$  with  $h$-map  $ {\mathbb Z}
\times   {\mathbb Z} \to   {\mathbb Z} $   given by   $(m,n) \mapsto
mn$, so that $h(1,1)$  represents the Hopf map  $\eta$.  But the
diagram (4.7) for the suspension  square of  $S^1$   may now be
completely determined, and is the universal  example for
Proposition 4.6.  This completes the proof of the proposition.
\end{proof}

 For the proof of the final part of Theorem
\ref{square} we have to explain how  the particular crossed module
given in the theorem determines the homotopy  2-type. This is proved
by considering the Moore complex of the diagonal simplicial group of
the bisimplicial group arising as the nerve of the  associated
cat$^2$-group. \nsec{Simply connected 3-types and crossed squares of
Abelian groups}     It is known that the 3-type of a simply
connected space  $X$  is  determined by the homotopy groups  $\pi_2
X ,~ \pi_3 X$  and the quadratic function  $\eta^{\ast}  : \pi_2 X
\to   \pi_3 X$  induced by composition with the Hopf  map $\eta  :
S^3 \to  S^2$.  This essentially results from the fact that for
abelian groups  $A$ and   $B$   the cohomology group
$H^4(K(A,2),B)$  is isomorphic to the group of quadratic  functions
$A \to   B$  (Eilenberg and Mac Lane, 1954).  The aim of this
section  is to show that all simply connected  3-types can be
modelled by a crossed square of Abelian groups.  It is not known if
simply connected $(n + 1)$-types can be modelled by crossed
$n$-cubes of Abelian groups.   \begin{thm}   Let  $C$  and  $D$  be
Abelian groups such that  $C$  is  finitely  generated,  and let  $t
:  C \to D$  be a quadratic function.  Then there is a crossed
square
$$ G \qquad 
\sqdiagram{L}{\lambda}{M}{\lambda '}{1}{M}{-1}{M} $$ of abelian
groups whose classifying space  $X = \B G$  satisfies $\pi_2  X
\cong C$, $ \pi_3 X \cong  D$  and such that these isomorphisms map
$\eta ^{\ast} $  to the  quadratic map $t$.
\end{thm}  \begin{proof}  The quadratic function  $t$  has first
to be extended to a  biadditive  map.  We use a slight
modification of a definition of Eilenberg and  Mac Lane, 1954,
$\S$ 18.

  Let  $t : C \to   D$  be a quadratic function on
Abelian groups  $C,D$.  A  {\em biadditive extension} of  $t$
is an abelian group  $M$  and an  epimorphism   $\alpha  : M \to
C$  of Abelian groups together with a biadditive map   $\phi  :
M \times  M \to   D$  such that for all  $m,m' \in  M$

(5.1.1)
$\phi (m,m) = t\alpha m$;

(5.1.2) $\phi (m,m') = 0$  if $\alpha
m = \alpha m' = 0$;

 (5.1.3) $\phi (m,m') = \phi (m',m)$.

\noindent It is
shown in {\em loc. cit.} that such a biadditive extension exists
assuming  $C$  is finitely generated.  (In fact they do not
assume the symmetry  condition (5.1.3), but their proof of
existence yields such a  $\phi$.)

  Let  $K = \Ker \alpha$
and let  $L$  be the product group  $D \times   K$.   Let  $M$ act
on  $L$  on the left by $$^m(d,k) = (d + \phi (m,k),k),$$ for  $m
\in  M , ~d \in  D , ~k \in  K$.  Define  $\lambda ,\lambda ' : L
\to    M$ by  $\lambda (d,k) = -k , \lambda '(d,k) = k$,  for $(d,k)
\in  L$, and let  $M$  act trivially on itself.  Then $\lambda$
and   $\lambda '$   are $M$-morphisms, and (5.1.2) shows that they
are also crossed modules.   Define  $h : M \times  M \to L$  by
$$h(m,m') = (\phi (m,m'),0)$$  for  $m,m' \in  M$.   A
straightforward check shows that we have defined a  crossed square
$G$  say.  The symmetry condition, or even the weaker  condition
that  $\phi (m,m') = \phi (m',m)$  if  $m$  or  $m'$  lies in   $K$,
is used to verify that  $$h(\lambda (d,k),m) = (d,k) -{} ^m(d,k).$$

The homotopy groups of  $\B G$   are computed as the homology groups
of the  chain complex  $$ L \llabto{2}{(-\lambda,\lambda ')} M
\times M \llabto{1}{\psi} M $$ where  $\psi (m,m') = m - m'$. Thus
$\pi_2 \B G \cong  M/K \cong  C ,~  \pi_3  \B G \cong  D$. Further
$h(m,m) = (\phi (m,m),0) = (t\alpha m,0)$.  This proves the final
assertion of the theorem.                \end{proof}

Note that by the proof of this theorem, while the groups are
Abelian,  the actions are in general non-trivial.  So the associated
cat$^2$-group in  general has non-Abelian big group.
\vspace{-4ex}\section*{Acknowledgments} \vspace{-2ex}    I would
like  to thank J.-L. Loday for conversations on the material of this
paper.  The work was supported by: the British Council; the
Universit\'e  Louis Pasteur, Strasbourg; and the SERC. \end{document}